\documentclass[11pt,spanish]{amsart} 
\usepackage{amssymb, amsthm, amsmath,mathrsfs} 
\usepackage{latexsym}
\usepackage[latin1]{inputenc}
\usepackage[T1]{fontenc}

\usepackage{xcolor}
\usepackage{tikz}
\usepackage{float}
\usepackage[pdftex]{hyperref}
\usepackage{enumerate}
\usepackage{soul}
\usepackage{comment}
\usepackage{booktabs}
\usepackage{orcidlink}

\hypersetup{
 colorlinks = true,
 urlcolor = blue!70!black,
 linkcolor = red!70!black,
 citecolor = green!70!black
}

\theoremstyle{plain}
\newtheorem{theorem}{Theorem}[section]
\newtheorem{lemma}[theorem]{Lemma}
\newtheorem{proposition}[theorem]{Proposition}

\theoremstyle{definition}

\theoremstyle{remark}
\newtheorem{remark}[theorem]{Remark}

\numberwithin{equation}{section}
\numberwithin{theorem}{section}

\title[Sobolev OP: Connection formulae]
{Sobolev orthogonal polynomials: Connection formulae}
\author[Roberto S. Costas-Santos]{Roberto S. Costas-Santos \orcidlink{0000-0002-9545-7411}}
\address[Roberto S. Costas-Santos]{\small \rm Dpto. de M\'etodos Cuantitativos, 
Universidad Loyola de Andaluc\'ia, 
E-41704, Dos Hermanas, Seville, Spain}
\email{rscosa@gmail.com}

\date{\today}


\begin{document}
\begin{abstract}
This contribution aims to obtain several connection formulae 
for the polynomial sequence, which is orthogonal with respect 
to the discrete Sobolev inner product
\[
\langle f, g\rangle_n=\langle {\bf u}, fg\rangle+ \sum_{j=1}^M \mu_{j} f^{(\nu_j)}(c_j) g^{(\nu_j)}(c_j),
\]
where ${\bf u}$ is a classical linear functional, $c_j\in \mathbb R$, $\nu_j\in \mathbb N_0$, 
$j=1, 2,...., M$. The Laguerre case will be considered.
\end{abstract}

\maketitle
\section{Introduction}
\noindent In this paper, we are going to consider sequences of polynomials orthogonal with 
respect to the discrete Sobolev inner product
\begin{equation} \label{eq:01}
\langle f, g\rangle=\langle {\bf u}, fg \rangle+ \sum_{j=1}^M \mu_{j} f^{(\nu_j)}(c_j) g^{(\nu_j)}(c_j),
\end{equation} 
where ${\bf u}$ is a classical linear functional, $c_j, \mu_{j}\in \mathbb C$, and $\nu_j\in \mathbb 
N_0$, $j=1, 2,...., M$ in the widest sense possible for all the parameters and values related to the 
polynomials we want to study. 

Observe that without loss of generality, we can assume that $\nu_1\le \nu_2\le \cdots 
\le \nu_M$. 
For this reason, we will call the polynomials orthogonal with respect to \eqref{eq:01} 
sequentially ordered Sobolev-type orthogonal polynomials. 

Observe one can express the inner product \eqref{eq:01} in the following compact 
way \cite{mr1320230}
\begin{equation} \label{eq:02}
\langle f, g\rangle=\langle {\bf u}, fg \rangle+({\mathbb D} f)^TD {\mathbb D} g,
\end{equation} 
where ${\mathbb D}$ is the vector differential operator defined as
\[
{\mathbb D} f:=\left(\left.f^{(\nu_1)}(x)\right|_{x=c_1},\left.f^{(\nu_2)}(x)\right|_{x=c_2}, 
..., \left.f^{(\nu_M)}(x)\right|_{x=c_M}\right)^T,
\]
$D$ is the diagonal matrix with entries $\mu_{1}, ..., \mu_{M}$ and $A^T$ is the transpose 
of the matrix $A$.

To simplify the notation we will write throughout the document $f^{(\nu)}(c)$ instead of 
$\left.f^{(\nu)}(x)\right|_{x=c}$.

\begin{remark}\label{rem:1}
Observe in the case when $\nu_1=0$, $\nu_2=1$, ..., $\nu_M=M-1$ and all the mass-points 
are equal each others the authors usually denote by $\mathbb F$ the matrix ${\mathbb D} f$ 
(see \cite{mr3197725} and references therein). 
\end{remark}
For a more detailed description of this Sobolev-type orthogonal polynomials (including the 
continuous ones) we refer the readers to the reviews \cite{mr1246854,mr3360352}.

The structure of the paper is as follows: in Section \ref{sec2}, some preliminary results are quoted. 
In section \ref{sec3}, all the algebraic results are presented, such as several connection formulas 
for the sequentially ordered balanced Sobolev-type orthogonal polynomials, a hypergeometric 
representation for these polynomials, as well as some other algebraic relations between the 
classical orthogonal polynomials and the discrete-Sobolev ones. 
\section{Auxiliary results}\label{sec2}
We adopt the following set 
notations: $\mathbb N_0:=\{0\}\cup\mathbb N=\{0, 1, 2, ...\}$, and we 
use the sets $\mathbb Z$, $\mathbb R$, $\mathbb C$ which represent 
the integers, real numbers, and complex numbers, respectively.
Let $\mathbb P$ be the linear space of polynomials and let $\mathbb P'$ be its algebraic dual space.

We will also adopt the following notation: We denote by $\langle {\bf u}, p\rangle$ the duality bracket 
for ${\bf u}\in \mathbb P'$ and $p\in \mathbb P$, and by $({\bf u})_n =\langle {\bf u}, x^n\rangle$, with 
$n\ge 0$, the canonical moments of ${\bf u}$.

For any $n\in \mathbb N_0$, $a\in \mathbb C$, the Pochhammer symbol, or shifted factorial, is 
defined as
\[
(a)_n:=a(a+1)\cdots (a+n-1).
\] 
The Taylor polynomial of degree $N$ is defined as 
\[
\left[f(x;c)\right]_N:=\sum_{k=0}^N \frac{f^{(k)}(c)}{k!}(x-c)^k,
\]
for every function $f$ for which $f^{(k)}(c)$, $k = 0,1,2, ...,N$ exists. 

The hypergeometric series is defined for $z\in \mathbb C$, $s, r \in \mathbb N_0$, 
$b_j\not \in -\mathbb N$ as \cite[\S 1.4]{mr2656096}
\[
{}_rF_s\left(\begin{array}{c} a_1, ..., a_s\\ b_1, ..., b_r\end{array} ; x \right)=
\sum_{k=0}^\infty \frac{(a_1)_k\cdots (a_s)_k}{(b_1)_k\cdots (b_r)_k}\frac{x^k}{k!}.
\]
Given a moment functional $\bf u$, it is  said to be quasi-definite or regular (see \cite{mr0481884}) if 
the Hankel matrix $H = \left(({\bf u})_{i+j}\right)_{i,j=0}^\infty$ associated with the moments of the 
functional is quasi-definite, i.e., all the $n$-by-$n$ leading principal submatrices are regular 
for all $n\in \mathbb N_0$. Hence, there exists a sequence of polynomials 
$\left(P_n\right)_{n\ge0}$ such that
\begin{enumerate} 
\item The degree of $P_n$ is $n$.
\item $\langle {\bf u}, P_n(x) P_m(x)\rangle =0$, $m\ne n$.
\item $\langle {\bf u}, P^2_n(x)\rangle =d^2_n\ne 0$, $n=0, 1, 2, ...$
\end{enumerate}
Special cases of quasi-definite linear functionals are the classical ones (Jacobi, Laguerre, 
Hermite and Bessel).

We denote the $n$-th reproducing kernel by
\[
K_n(x,y)=\sum_{k=0}^n \frac{P_k(x){P_k(y)}}{d_k^2}.
\]
From the Christoffel-Darboux formula (see \cite{mr0481884} or \cite[Eq. (3.1)]{MR2500498}), we have
\begin{equation} \label{eq:03}
K_n(x,y)=\frac{k_n}{k_{n+1}d^2_n} \frac{P_{n+1}(x){P_n(y)}-P_n(x){P_{n+1}(y)}}{x- y},
\end{equation} 
where $k_m$ is leading coefficient of $P_m(x)$, $m\in \mathbb N_0$.

We will use the following notation for the partial derivatives of $K_n(x, y)$:
\[
\frac{\partial^{j+k}}{\partial x^j\partial y^k}K_n(x,y)=K^{(j,k)}_n(x,y),\quad 0\le j, k\le n.
\]
Note that, when $j = k = 0$, $K_{n}(x, y)=K^{(0,0)}(x, y)$ is the usual 
reproducing Kernel polynomial.

A direct consequence of the Christoffel-Darboux formula \eqref{eq:03} is 
the following result \cite[Eq. (16)]{MR2581392}:
\begin{proposition} \label{pro:1}
The $j$-th partial derivative of the $n$-th reproducing kernel can be written as
\begin{equation} \label{eq:04}
K^{(j,0)}_n(x,y)=\frac{k_n\, j!}{k_{n+1}d_n^2}\frac{P_{n+1}(x) {[P_{n}(x;y)]_j}-P_n(x) {[P_{n+1}(x;y)]_j}}
{(x-\overline{y})^j}.
\end{equation}
\end{proposition}  
Observe the following consequence, provided that $c$ is not a zero of $P_n(x)$
for any $n$
\begin{equation} \label{eq:05}
\frac{\left|P^{(j)}_n(c)\right|^2}{d^2_n}=K^{(j,0)}_n(c,x)-K^{(j,0)}_{n-1}(c,x)
\end{equation}
One last result will be useful to obtain some of our algebraic results.
\begin{lemma} \label{lem:1} \cite[Lemma 2.1]{mr1391618}
Let $M\in \mathbb N$, ${\bf u}$ be a classical linear form. Let $c_1, c_2, ..., c_M\in \mathbb R$, 
$\nu_1, \nu_2, ...,\nu_M\in \mathbb N_0$, and let us denote by $(S^{\vec \mu}_n(x;\vec \nu,\vec c))$ 
the sequence of orthogonal polynomials with respect to the inner product \eqref{eq:01}. 
If $c_i$ is not a zero of $S^{\vec \mu}_n(x;\vec \nu,\vec c)$, $i=1, 2, ..., M$ for all $n\in \mathbb N_0$ 
then, there exists a polynomial, namely $\zeta(x)$, such that 
${\mathbb D} \left(\zeta(x)S^{\vec \mu}_n(x;\vec \nu,\vec c)\right)=\vec 0$ holds. 
\end{lemma} 
\begin{remark}\label{rem:2} Observe that if all the $c_i$'s are all different then 
$\zeta(x)=\prod_{j=1}^M (x-c_j)^{\nu_j+1}$, and if all of them are equal to each other, i.e.  $c_i=c$ for 
$i=1, 2, ..., M$, then $\zeta(x)=(x-c)^{\nu_M+1}$.
\end{remark}
Without loss of generality, we denote by $\zeta (x)$ to the polynomial of minimum degree among all nonzero 
polynomials satisfying the conditions of the Lemma \ref{lem:1}.
\subsection{The Laguerre polynomials}
Let $\left(L_n^{\alpha}(x)\right)$ be the sequence of Laguerre polynomials, orthogonal 
with respect to the linear form ${\bf u}_\alpha$ on $\mathbb P$. These polynomial sequence 
is classical since ${\bf u}_\alpha$ fulfills the Pearson equation
\[
\frac d{dx}\left[x \, {\bf x}\right]=(\alpha+1-x)\,{\bf x}.
\]
\begin{remark}\label{rem:3}
Note that if $\Re(\alpha)>-1$ then, the linear from ${\bf u}_\alpha$ has 
the following integral re\-pre\-sen\-tation (see for instance \cite{mr1405981}, \cite[\S 18.3]{dlmf} or  \cite[\S 9.12]{mr2656096}):
\[
\langle {\bf u}_\alpha, f\rangle=\int_0^\infty f(x)\, x^\alpha e^{-x} dx,
\]
and when $\alpha<0$ the orthogonality of Sobolev-type is given in \cite{mr1405981}.
\end{remark}
The Laguerre polynomial can be explicitly given in terms of hypergeometric series as 
\begin{equation} \label{eq:06}
L_n^{\alpha}(x)=\frac{(\alpha+1)_n}{n!}{}_1F_1\left(\begin{array}{c}-n \\ \alpha+1\end{array} ; x\right).
\end{equation}

Let us summarize some basic properties of the Laguerre orthogonal polynomials that will be used throughout 
this work. 
\begin{proposition} \label{prop:2} Let $(L_n^{\alpha}(x))$ be the  Laguerre polynomials. 
The following statements hold:
\begin{itemize} 
\item The three-term recurrence relation:
\begin{equation} \label{eq:07}
(n+1)L_{n+1}^{\alpha}(x)+(x-2n-\alpha-1)L_n^{\alpha}(x)+(n+\alpha)L_{n-1}^{\alpha}(x),\quad 
n=1, 2, ...,
\end{equation}
with initial conditions $L_0^{\alpha}(x)=1$ and $L_1^{\alpha}(x)=\alpha+1-x$.
\item The first structure relation:
\begin{equation} \label{eq:08}
x \left(L_n^{\alpha}(x)\right)'=n L_n^{\alpha}(x)-(n+\alpha) L_{n-1}^{\alpha}(x),\quad 
n=1, 2, ...
\end{equation}
\item The second structure relation:
\begin{equation} \label{eq:09}
L_n^{\alpha}(x)=-\left(L_{n+1}^{\alpha}(x)\right)'+\left(L_{n}^{\alpha}(x)\right)', \quad 
n=0, 1, ...
\end{equation}
\item The squared norm:
\begin{equation} \label{eq:10}
d^2_n=\frac{(\alpha+1)_n}{n!},\quad 
n=0, 1, ...
\end{equation}
\item The Ladder operators:
\begin{eqnarray} 
\left(L_{n}^{\alpha}(x)\right)'&=&-L_{n-1}^{\alpha+1}(x), \quad n=1, 2, ...,\label{eq:11}\\
x \left(L_{n}^{\alpha}(x)\right)'+(\alpha-x) L_{n}^{\alpha}(x)&=&(n+1)L_{n+1}^{\alpha-1}(x), \quad n=0, 1, ...\label{eq:12}
\end{eqnarray}
\end{itemize}
\end{proposition}

\section{The connection formulae} \label{sec3}
A first step to get asymptotic properties is to obtain an adequate expression of the 
polynomials $L_n^{\alpha, \vec \mu}(x)$ in terms of the Laguerre 
polynomials, i.e., to solve the connection problem. 

\begin{remark}\label{rem:4}
\begin{itemize} 
\item[] 
\item Observe that, by construction, it is clear that $L_n^{\alpha, \vec \mu}(x)=L_n^{\alpha}(x)$
for $n=0, 1, ..., \nu_1-1$.
\item For the part of algebraic calculations in this work, the assumption of the dependency of the parameters 
$\mu_j$ of $n$ is unnecessary. Therefore we will omit such dependency.
\item We assume that $L_n^{\alpha, \vec \mu}(x)$ has the same leading coefficient than $L_n^{\alpha}(x)$.
\end{itemize} 
\end{remark}
Since the Laguerre polynomials constitute a basis of the polynomials,  
we can consider the Fourier expansion of $L_n^{\alpha, \vec \mu}(x)$ in terms of such polynomial sequence.
\begin{proposition}\label{prop:3}
For every $n\ge \nu_1$ the following identity holds:
\begin{equation} \label{eq:14}
L_n^{\alpha, \vec \mu}(x)=L_n^\alpha(x)-\sum_{j=1}^M \mu_{j}
\left(L_n^{\alpha, \vec \mu}(c_j)\right)^{(\nu_j)}K^{(0,\nu_j)}_{n-1}(x,c_j).
\end{equation}  
\end{proposition}
This is a classical result, so the proof will be omitted. Still, we need to compute the values 
$\left(L_n^{\alpha, \vec \mu}(c_j)\right)^{(\nu_j)}$ for $j=1, 2, ..., M$ in order to have the complete expression. 
We must use an analogous result from \cite[Proposition 2]{mr1357591}. If we define
$\mathbb L_n$ as ${\mathbb D}L_n^{\alpha}(x)$,  $\mathbb S_n$ as ${\mathbb D}L_n^{\alpha,\vec \mu}(x)$ and
$\mathbb K_n={\mathbb D}^T_x{\mathbb D}_y K_n(x,y)$, then 
we need to solve the linear system
\[
\mathbb S_n=\mathbb L_n-{\mathbb K}^T_{n-1} D^T \mathbb S_n.
\]
Therefore, after some straightforward manipulations, \eqref{eq:14} becomes
the desired compact connection formula \cite[Proposition 2]{mr1357591}.
\begin{equation} \label{eq:15} 
L_n^{\alpha, \vec \mu}(x)=L_n^{\alpha}(x)-\mathbb L^T_n(\mathbb I
+D {\mathbb K}_{n-1})^{-1}D {\mathbb K}_{n-1}(x),
\end{equation} 
where ${\mathbb K}_{n}(x)={\mathbb D}_y K_n(x,y)$.

We expect to obtain this identity, but there are other connection formulas.
\begin{remark}\label{rem:5}
Observe that the discrete Laguerre-Sobolev polynomials exits for all $n$ if and only if the matrices
 $\mathbb I+D {\mathbb K}_{n-1}$ are regular for all $n=1, 2, ...$.
\end{remark}
In the next result, we establish a connection formula for the discrete Sobolev polynomials 
$L_n^{\alpha, \vec \mu}(x)$ similar to the one obtained in \cite[Theorem 1]{mr3345273} for  non-varying 
discrete Sobolev orthogonal polynomials. 

Let $\zeta(x)$ be the polynomial, of degree $\nu$, we obtain from Lemma \ref{lem:1}, then 
it is clear that for any two polynomials $f$ and $g$, we have 
\begin{equation} \label{eq:16}
\langle \zeta(x) f(x), g(x)\rangle=\langle {\bf u}_\alpha, \zeta(x) f(x) \overline{g(x)}\rangle=\langle f(x), \overline{\zeta(x)} g(x)\rangle.
\end{equation} 

\begin{proposition}
Let  $(\zeta_j(x))_{j=0}^\nu$ be a sequence of polynomials, with $\zeta_\nu(x)=\zeta(x)$, such that $\deg 
\zeta_k(x)=k$ and it is a divisor of $\zeta_{k+1}(x)$ for $k=0, 1, ..., \nu-1$, and let $\big(P_n^{[\zeta^2_j]}(x)\big)$ 
be the polynomials orthogonal with respect to the linear functional $|\zeta_j(x)|^2{\bf u}_\alpha$, for $j=0, 1, ..., \nu$.

If the following conditions hold
\begin{equation} \label{eq:17}
P_n(c_j) P_{n-1}^{[\zeta^2_1]}(c_j)\cdots P_{n-\nu}^{[\zeta^2_\nu]}(c_j)\ne 0, \quad j=1, 2, ...., M,
\end{equation} 
then, there exists a family of coefficients $(\lambda_{j,n})_{j=0}^{\nu}$, not identically zero, such that 
for any $n\ge \nu$  the following connection formula holds:
\begin{equation} \label{eq:18}
L_n^{\alpha, \vec \mu}(x)=\sum_{j=0}^{\nu} \lambda_{j,n}\, \zeta_j(x) P_{n-j}^{[\zeta^2_j]}(x).
\end{equation}  
\end{proposition}
Another connection formula connects the discrete Sobolev polynomials with the 
the derivatives of the Laguerre polynomials, which proof is similar to the previous one.
\begin{proposition}
For every $n\ge \nu$ the following identity holds:
\begin{equation} \label{eq:20}
L_n^{\alpha, \vec \mu}(x)=\sum_{k=0}^{\nu} \xi_{k,n}\, L^{\alpha+k}_{n-k}(x).
\end{equation}  
\end{proposition}

Extending these results to other classical families is straightforward, even more 
generic frameworks such as discrete classical polynomials.
\bibliographystyle{plainurl}

\end{document}